\documentclass[12pt]{amsart}
\usepackage{latexsym, amsmath, amscd, amssymb, amsthm} 
\usepackage[T2A]{fontenc}
\usepackage[cp1251]{inputenc}
\newtheorem{te}{Теорема}

 \newtheorem*{oz}{Означення}
 
 \newtheorem*{lm}{Лема}

\begin{document}

 \title{ Симетричнi iнварiанти деяких модулярних алгебр Лi картанiвського типу }

\author{Л. Бедратюк}\address{ Хмельницький національний університет, вул. Інститутська, 11, 29016, м. Хмельницький}

\begin{abstract}
 Let  $L$   be one  of the finite dimensional  Lie algebras  $W_n({\bf m}),$ $S_n({\bf m}),$ $  H_n({\bf m})$  of  Cartan type over an algebraically closed field of prime characteristic  $p>0.$     For  an elements $F$ of the symmetrcial algebra $S(L)$  we   found necessary and sufficient  condition  in order to  the element  $ad(\partial_1)^{p^{m_1}-1} ad(\partial_2)^{p^{m_2}-1}\cdots ad(\partial_n)^{p^{m_n}-1}(F)$  belongs to the symmetrical invariants algebra $S(L)^L.$ Also,  for $p=3,5$     the algebra of symmetrical invariants $S(H_2)^{H_2}$  is  calculated  in  explicit   way. 
\end{abstract}
\maketitle

{\bf 1.} Нехай $L$ одна із простих скінченовимірних алгебр Лі картанівського типу $W_n({\bf m}),$ $S_n({\bf m}),$ $ H_n({\bf m}),$ де  ${{\bf m}=(m_1,m_2,\ldots, m_n),} $ $m_i \in \mathbb{Z},$ яка розглядається   над алгебраїчно замкнутим полем $\mathbb{K}$ додатньої характеристики $p.$  Алгебра $L$  є градуйована 
$$
L=L_{-1} \oplus L_0 \oplus \cdots \oplus L_r, [L_i,L_j] \subset L_{i+j}.
$$
Тут $r=n \sum_k(p^{m_i}-1)-1$  для  $W_n({\bf m})$ і  $r=n \sum_k(p^{m_i}-1)-2$  для алгебр $S_n({\bf m}),$ $  H_n({\bf m}).$
Компоненти відповідної фільтрації позначимо через $\mathcal{L}_{i}=L_i \oplus L_{i+1} \oplus \cdots \oplus L_r, $  $\mathcal{L}_{-1}=L,$ див. всі деталі в \cite{KS}.
Обчислення центру $Z(L)$ універсальної  огортуючої алгебри $U(L)$ є однією із важливих задач у  відкритій проблемі класифікації  незвідних зображень алгебри $L.$ На сьогоднішній день центр $Z(L)$ повністю описаний лише для алгебри  $W_1({\bf m})$  в \cite{Er}. Окремi центральнi елементи для  деяких алгебр знайдені  в  \cite{Kor}, \cite{Dz}.

Першим кроком в описі  алгебри $Z(L)$ може  бути знаходження алгебри  інваріантів $S(L)^L$   симетричної алгебри $S(L)$ відносно приєднаної дії алгебри $L.$  Алгебра  $S(L)^L$  є скінченно породженою алгеброю над $p$-центром алгебри $L.$  Мімальна кількість породжуючих елементів алгебри $S(L)^L$  обчислена в роботі \cite{Kr}, проте явний вигляд цих елементів в загальному випадку невідомий. 

В  роботах    \cite{Bk-1},  \cite{Bk-2}  дано опис алгебри  $S(L)^{L_{-1}}.$  Зокрема, було показано, що всякий нетривіальний однорідний симетричний інваріант  $S(L)^L,$степінь якого не ділиться на характеристику поля,  записується у вигляді $d^{(\delta)}(F),$ де $$d^{(\delta)}:=ad(\partial_1)^{p^{m_1}-1} ad(\partial_2)^{p^{m_2}-1}\cdots ad(\partial_n)^{p^{m_n}-1},$$  a $F \in S(L),$ $\partial_1, \partial_2,\ldots, \partial_n$ -- базис простору $L_{-1}.$

В  даній  роботі встановлено  умови яким повинен задовольняти елемент $F \in S(L)$  для того щоб елемент  $d^{(\delta)}(F)$  був  симетричним інваріантом алгебри  $L.$  Показано,  що для алгебр $S_n({\bf m}),$ $ H_n({\bf m}),$ $ \overline H_n({\bf m}),$ ($W_n({\bf m})$)  елемент $F$ має бути інваріантом (відповідно напівінваріантом)  максимальної підалгебри $\mathcal{L}_{0}.$ Використовуючи отриманий критерій, обчислено алгебру симетричних інваріантів для $H_2:=H_{2}({\bf m}),$  при $ {\bf m}=(1,1,\ldots,1)$  і $p=3,5.$

${\bf 2.}$ Нагадаємо  означення алгебр Лі  $W_n({\bf m}), S_n({\bf m}), H_n({\bf m}).$  Зафіксуємо позначення  для наборів з $\mathbb{Z}_{+}^n$ --
 $\delta:=(\delta_1, \delta_2,\ldots,\delta_n), \delta_i=p^{m_i-1},$  $\epsilon_i :=(0,\ldots,1,\ldots,0).$ Покладемо  $|\alpha|:=\sum_i \alpha_i,$ $\alpha \in \mathbb{Z}_{+}^n.$

Алгеброю розділених степенів $\mathbb{K}_n({\bf m})$ називається комутативна алгебра, яка задана твір\-ними  $x_1,$ $x_2, \ldots,$ $x_n,$  співвідношеннями $x_1^{p^{m_1}}=0,$ $\ldots $ $x_n^{p^{m_n}}=0$  і правилом множення
$$
x^{(\alpha)} x^{(\beta)} ={\alpha+\beta \choose \alpha} x^{(\alpha+\beta)}, {\alpha \choose \beta}=\prod_i  {\alpha_i \choose \beta_i}
$$
тут  $x^{(\alpha)}:=x_1^{\alpha_1} \ldots, x_n^{\alpha_n},$  $\alpha \in \mathbb{Z}_{+}^n,$ $\alpha \leq \delta.$

Загальна алгебра $W_n({\bf m})$  породжується  всіма спеціальними диференціюваннями алгебри $\mathbb{K}_n({\bf m})$  вигляду 
$$
D=\sum_i f_i \partial_i, f_i \in k_n({\bf m}).
$$
Спеціальна алгебра $S_n({\bf m}), n \geq 2$ є підалгеброю алгебри $W_n({\bf m}),$  яка  породжена диференціюваннями 
$$
\mathcal{D}_{i,j}(\alpha)=\partial_i(x^{(\alpha)})\partial_j -\partial_j(x^{(\alpha)})\partial_i, i<j \leq n, \alpha \in  \mathbb{Z}_{+}^n.
$$
Гамільтонова алгебра $H_n({\bf m}),$ $n$-парне, складається із диференціювань 
$$
D(\alpha)=\sum_{i=1}^n a_{i,\pi i} \partial_i(x^{(\alpha)})\partial_{\pi i}, \alpha \in \mathbb{Z}_{+}^n,  \alpha < \delta,
$$
де $\pi$  --  інволютивна перестановка без нерухомих точок множини $\{ 1,2,\ldots,n \},$  причому  ${a_{i,\pi i}=\pm 1}$, ${a_{i,\pi i}+a_{\pi i, i}=0}.$  Всі вищеозначені алгебри є простими алгебрами Лі.
Крім вказаних алгебр  ми  також будемо  працювати з алгеброю  $\overline H_n({\bf m}):=H_n({\bf m}) \oplus \langle D(\delta) \rangle.$
Алгебра  симетричних інваріантів $S(L)^L$ визначається як анулятор $L$-модуля $S(L).$

${\bf 3.}$ Для кожної алгебри $W_n({\bf m}), S_n({\bf m}), H_n({\bf m}),\bar H_n({\bf m}) $ компонента $L_{-1}$ породжена комутуючими диференціюваннями $\partial_1,$ $\ldots,$  $\partial_n,$ тому оператор $d^{(\delta)}$  визначений коректно для всіх типів алгебр.
\begin{oz} 
 Елемент $F\in S(L)$  називається {\bf генератором } симетричного інваріанта $z,$ якщо виконується рівність  $z=d^{(\delta)}(F).$
\end{oz} 
Знання  генератора  симетричного   інваріанта важливе з обчислювальної  точки зору, оскільки, генератор  має   набагато  простіший вигляд  ніж відповідний  йому симетричний  інваріант. 
В  роботі \cite{Bk-1} встановлено,  що для кожного  однорідного симетричного інваріанта $z,$ алгебри Лі $L,$ для якого ${\deg(z) \not= 0  \mod p},$
існує єдиний, з точністю до ядра оператора $d^{(\delta)},$ генератор.
 Наступна  теорема встановлює якими властивостями  володіє генератор симетричного інваріанта. 

\begin{te} 
Нехай $z$  -- нетривіальний симетричний однорідний інваріант алгебри $L,$  $  \deg(z) \not=  0 \mod p,$ і  $F$ -- його генератор. Тоді
\begin{itemize}
\item[ $(i)$]  якщо  $L=W_n({\bf m})$   то $ \mathcal{L}_{1}(F)=0$  i $ad(-x_i \partial_j)(F) = \delta_{i,j} F ,$ $\delta_{i,j}$  --символ Кронекера;
\item[ $(ii)$]  якщо $L= S_n({\bf m}) $ або   $L=\overline  H_n({\bf m}),$ то $ \mathcal{L}_{0}(F)=0.$
\end{itemize}
\end{te} 
\begin{proof}[Доведення]
$(i)$  В роботі \cite{Bk-2}   показано, що елемент  $F$  можна  подати у вигляді 
$$
F=z_1A_1+z_2 A_2+\cdots +z_n A_n,
$$
де $z_i=x^{(\delta)} \partial_i,$ $A_i \in S(L)$, а векторний простір $A$, натягнутий на елементи $A_1, A_2, \ldots A_n$  є $ \mathcal{L}_{1}$-модулем з нульовою дією. Крім того $A$ є   незвідним $L_0$-модулем  з наступною дією 
$$
D_{i,j}(A_k)=-\delta_{j,k} A_i, D_{i,j}:=-ad(x_i \partial_j).
$$
Враховуючи те, що $D_{i,j}(z_k)=\delta_{i,j} z_k +\delta_{i,k} z_j$ отримаємо 
$$
\begin{array}{l}
\displaystyle D_{i,j}(F)= D_{i,j}\Bigl(\sum_k z_k A_k\Bigr)=\sum_k\Bigl(D_{i,j}(z_k) A_k+z_k  D_{i,j}(A_k)\Bigr)=\\

\displaystyle=\sum_k\Bigl((\delta_{i,j} z_k {+}\delta_{i,k} z_j A_k{+}z_k (-\delta_{j,k} A_i)\Bigr){=} \delta_{i,j} F{+} z_j  \sum_k \delta_{i,k} A_k {-}A_i \sum_k \delta_{j,k} z_k =
\\
\displaystyle= \delta_{i,j} F+z_j A_i-z_j A_i= \delta_{i,j} F.
\end{array}
$$
Оскільки, для всіх $ i$ маємо $\mathcal{L}_{1}(z_i)=0$ і $\mathcal{L}_{1}(A)=0,$ то і $\mathcal{L}_{1}(F)=0.$

$(ii)$
Нехай $L=S_n({\bf m}).$
В  роботі \cite{Bk-1} показано, що довільний  однорідний симетричний інваріант  ${z \in S(L)^L,}$ $ \deg(z)\not= 0 \mod p$  записується у спеціальному вигляді $z=d^{(\delta)}(F)$, причому $F=\sum_{i<j} z_{i,j}(\delta) A_{i,j}$,
де $A_{i,j} \in S(L)$    і $z_{i,j}:=\partial_i(x^{(\delta)})\partial_j-\partial_j(x^{(\delta)})\partial_i.$
  Покладемо $D_{i,j}=ad(x_i \partial_j), $ $T_{i,j}=ad(x_i \partial_j-x_j\partial_j)$. Як і в роботі \cite{Bk-2} можна довести наступне твердження
\begin{lm} 
Векторний простір $A$  породжений елементами $A_{i,j}$  є  $\mathcal{L}_{0}$-модулем з наступною дією 
$$
D_{s,t}(A_{i,j})=\delta_{t,j} A_{i,s}+\delta_{i,t} A_{s,t}, T_{s,t}(A_{i,j})=(\delta_{j,s}+\delta_{i,s}-\delta_{i,t}-\delta_{j,t}) A_{i,j},
$$
i $\mathcal{L}_{1}(A)=0$
\end{lm}

Прямими обчисленнями в алгебрi  $L$  отримуємо:
$$
\begin{array}{l}
D_{s,t}(z_{i,j})=[x_s \partial_t,\partial_i(x^{(\delta)})\partial_j-\partial_j(x^{(\delta)})\partial_i]=\delta_{i,s} z_{j,t}+\delta_{j,s} z_{t,i},\\
T_{s,t}(z_{i,j})=[x_s \partial_s-x_t \partial_t,\partial_i(x^{(\delta)})\partial_j-\partial_j(x^{(\delta)})\partial_i]=(\delta_{i,s}+\delta_{j,s}-\delta_{i,t}-\delta_{j,t}) z_{j,i}.
\end{array}
$$

Знаходимо: 
$$
\begin{array}{l}
\displaystyle D_{s,t}(F)=D_{s,t}\Bigl( \sum_{i<j} z_{i,j} A_{i,j}\Bigr)=\sum_{i<j}\Bigl(\bigl(\delta_{i,s} z_{j,t}+\delta_{j,s} z_{t,i}A_{i,j}+z_{i,j} (\delta_{t,j} A_{i,s}+\delta_{i,t} A_{s,t})\bigr)\Bigr)=\\
\displaystyle=\sum_{i<j}\delta_{i,s} z_{j,t} A_{i,\,j}+\sum_{i<j}\delta_{j,s} z_{t,i}A_{i,j}+\sum_{i<j}\delta_{t,j}z_{i,j}  A_{i,s}+\sum_{i<j}\delta_{i,t} A_{s,t}z_{i,j}=\\
\displaystyle=\sum_{s<j} z_{j,t} A_{s,\,j}+\sum_{i<s} z_{t,i}A_{i,s}+\sum_{i<t}z_{i,t}  A_{i,s}+\sum_{t<j} A_{s,j}z_{t,j}=\\
\displaystyle=\sum_{s<i} z_{i,t} A_{s,i}+\sum_{i<s} z_{t,\,i}A_{i,s}+\sum_{i<t}z_{i,t}  A_{i,\,s}+\sum_{t<i} A_{s\,,i}z_{t,i}=\\
\displaystyle=\sum_{i} z_{t,i} A_{i,s}+\sum_{i} z_{i,t} A_{i,s}=\sum_{i} (z_{t,i}+z_{i,t}) A_{i,s}=0.
\end{array}
$$
і
$$
\begin{array}{l}
\displaystyle T_{s,t}\Bigl( \sum_{i<j} z_{i,j} A_{i,j}\Bigr){=}\sum_{i<j}\bigl((\delta_{i,s}{+}\delta_{j,s}{-}\delta_{i,t}{-}\delta_{j,t}) z_{j,i} A_{i,j}{+}z_{i,j} (\delta_{j,s}{+}\delta_{i,s}{-}\delta_{i,t}{-}\delta_{j,t}) A_{i,j}\bigr){=}\\
\displaystyle=\sum_{i<j} (\delta_{j,s}+\delta_{i,s}-\delta_{i,t}-\delta_{j,t})(z_{j,i}A_{i,j}+z_{i,j}A_{i,j})=0.
\end{array}
$$
Оскільки,  диференціювання $x_s \partial_t$ i $ x_s \partial_s-x_t \partial_t$  породжують алгебру $L_0$  то $L_0(F)=0.$  Врахувавши $\mathcal{L}_{1}(A)=0$  i $\mathcal{L}_{1}(z_{i,j}(\delta))=0$ отримаємо,  що $\mathcal{L}_{1}(F)=0.$ Тому і $\mathcal{L}_{0}(F)=0.$

$(iii)$  Нехай $L=\overline H_n({\bf m}).$   З роботи \cite{Bk-1}  випливає що $F \in \overline{H}_n({\bf m}),$ причому $F= u A,$ $A \in S( \overline{H}_n({\bf m})),$  тут $u=x^{(\delta)}.$  Тоді,  врахувавши ${\delta \choose \alpha} =(-1)^{|\alpha|}  \mod p,$ отримаємо
$$
 z:=d^{(\delta)}(F)=d^{(\delta)}(u A)=\sum_{\alpha} (-1)^{|\alpha|} d^{(\alpha)}(u) d^{(\delta-\alpha)}( A).
$$
Оскільки елементи $\{  d^{(\alpha)}(u), 0\leq \alpha < \delta, \}$ очевидно, утворюють базис алгебри ${H}_n({\bf m}),$ то елементи $\{ d^{(\delta-\alpha)}( A),  0\leq \alpha < \delta \}$, утворюють базис спряженого модуля ${H}_n({\bf m})^*.$ Але для гамільто\-нової алгебри існує ізоморфізм ${H}_n({\bf m})$-модулів $\varphi :   {H}_n({\bf m}) \to {H}_n({\bf m})^*,$ який  продовжується  до гомоморфізму  $\varphi :   {\overline H}_n({\bf m}) \to {\overline H}_n({\bf m})^*,$ причому $\varphi(u)=A.$  Оскільки $\mathcal{L}_{0}(u)=0$ то і $\mathcal{L}_{0}(A)=0$ звідки отримаємо $\mathcal{L}_{0}(F)=0.$

\end{proof}

Перехід  від знаходження   симетричного інваріанта   до  знаходження  його генератора  є важливим з обчислювальної точки зору, оскільки симетричні інваріанти в розгорнутому вигляді задаються дуже громіздкими виразами.

Доведемо обернену теорему.

\begin{te} 
Нехай $F \in S(L)$ причому:
\begin{itemize}
\item[ $(i)$]  для   $L=W_n({\bf m})$   виконується  $ \mathcal{L}_{1}(F)=0$  i $ad(x_i \partial_j)(F) = - \delta_{i,j} F, $
\item[ $(ii)$]  для  $L=S_n({\bf m}), $ або   $L=\overline H_{n}({\bf m}),$ виконується $ \mathcal{L}_{0}(F)=0.$
\end{itemize}
Тоді елемент $d^{(\delta)}(F)$ є   симетричнм інваріантом.
\end{te} 
\begin{proof}[Доведення]
Алгебра Лі $L$  породжується як алгебра диференціюваннями з $L_{-1}$ і $L_r.$  Диференціювання    $d_i=ad(\partial_i), \partial_i  \in L_{-1}$  комутують між собою, і,  врахувавши ${d_i^{p^{m_i}}=0}$, отримаємо, що   $d_i ( d^{(\delta)}(F))=0.$ Отже, для перевірки того, що елемент $d^{(\delta)}(F)$ є симетричним інваріантом достатньо показати, що ${L_r( d^{(\delta)}(F))=0.}$

Нехай $L=W_n({\bf m}).$   Компонента $L_r$ породжена диференціюваннями вигляду $x^{(\delta)} \partial_i.$ Покладемо $D_i=ad(x^{(\delta)} \partial_i)$  i покажемо, що  в умовах теореми  $D_i(d^{(\delta)}(F))=0$ для всіх $i \leq n.$ Маємо, 

$$
D_i \bigl( d^{(\delta)}(F)\bigr)=\sum_{\gamma} (-1)^{|\gamma|} {\delta \choose \gamma}  d^{(\delta-\gamma)}\bigl( d^{(\gamma)}(D_i)(F)\bigr).
$$
Для тих значеннь $\gamma$  для яких  $d^{(\gamma)}(D_i) \in  \mathcal{L}_{1}$ відповідні доданки суми рівні нулю.  Тому
$$  
\begin{array}{l}
\displaystyle D_i \bigl( d^{(\delta)}(F)\bigr)=(-1)^{|\delta|} {\delta \choose \delta} d^{(\epsilon_k)}(F)+\sum_j (-1)^{|\delta-\epsilon_j|} {\delta \choose \delta-\epsilon_j}  d^{(\delta-\epsilon_k)}(D_i)(F)=\\

\displaystyle =(-1)^n d^{(\epsilon_k)}(F)+\sum_j (-1)^{n-1} (-1)^{|\epsilon_j|} d^{(\epsilon_j)}((ad(x^{(\epsilon_j)}\partial_k)(F))=\\

\displaystyle =(-1)^n d^{(\epsilon_k)}(F)+(-1)^n d^{(\epsilon_k)}(ad(x^{(\epsilon_k)}\partial_k)(F))=d^{(\epsilon_k)}(F)((-1)^n+(-1)^{n+1})=0.
\end{array}
$$
Отже, $d^{(\delta)}(F)  \in S(L)^L.$

Нехай $L=S_n({\bf m}).$  У  цьому випадку компонента $L_r$  породжена всіма диференціюваннями вигляду ${\partial_i(x^{(\delta)}) \partial_j  -\partial_j(x^{(\delta)}) \partial_i.}$ Покажемо, що коли $F$ задовольняє  умовам теореми, то     ${D_{i,j}(\delta)\bigl(d^{(\delta)}(F)\bigr)=0},$ для всіх $i <j \leq n ,$  ${D_{i,j}(\delta)=ad(\partial_i(x^{(\delta)}) \partial_j  -\partial_j(x^{(\delta)}) \partial_i)).}$
Маємо 
$$
D_{i,j}(\delta) \bigl( d^{(\delta)}(F)\bigr)=\sum_{\gamma} (-1)^{|\gamma|} {\delta \choose \gamma}  d^{(\delta-\gamma)}\bigl( d^{(\gamma)}(D_{i,j}(\delta))(F))\bigr).
$$
Якщо, $|\gamma|>1,$ то за умовою $d^{(\gamma)}(D_{i,j})(F) =0.$ Тому ненульові доданки будуть лише  при таких наборах -- $\delta, \delta{-}\epsilon_1,\ldots,\delta{-}\epsilon_n .$  Але $d^{(\delta)}(D_{i,j}(\delta)) =D_{i,j}(0)=0.$ Тому
$$
D_{i,j}({\delta}) \bigl( d^{(\delta)}(F)\bigr)=\sum_k (-1)^{(\delta-\epsilon_k)}d^{\epsilon_k}\bigl(d^{(\delta-\epsilon_k)}(\delta)(F)\bigr)=\sum_k (-1)^{n-1} d^{{\phantom{}}^{(\epsilon_k)}}D_{i,j}(\epsilon_k)(F)=
$$
$$
=(-1)^{n-1} \sum_k d^{(\epsilon_k)}\bigl( ad(\partial_i(x^{\epsilon_k})\partial_j-\partial_j(x^{\epsilon_k})\partial_i)(F)\bigr)=(-1)^{n-1} \sum_k d^{(\epsilon_k)}(ad(\delta_{i,k} \partial_j-\delta_{j,k} \partial_i)(F))=
$$
$$
=(-1)^{n-1}\bigl(d^{(\epsilon_i)}(\partial_j(F))-d^{(\epsilon_j)}(\partial_i(F))\bigr)=(-1)^{n-1}\bigl(d^{(\epsilon_i)}d^{(\epsilon_j)}(F)-d^{(\epsilon_j)}d^{(\epsilon_i)}(F)\bigr)=0. 
$$
Таким чином $d^{(\delta)}(F)  \in S(L)^L.$

Для алгебри $L=\overline H_{n}({\bf m})$ доведення аналогічне.
\end{proof}
Зауважимо, що із даного доведення зовсiм не випливає  нетривіальність симетричного інваріанта $d^{(\delta)}(F),$ навіть  якщо $F$ задовольняє умовам теореми. Тому при   обчисленнях потрібно  окремо перевіряти чи  елемент $d^{(\delta)}(F)$  відмінний від нуля. 

Таким  чином,  ми звели  проблему  обчислення симетричних інваріантів до еквівалентної задачі обчислення  генераторів симетричних інваріантів, яка є простішою з  обчислювальної  точки зору.

${\bf 4.}$ Використаємо  {теорему 2} для знаходження серії  нетривіальних симетричних інваріантів алгебр $\overline H_n({\bf m})$ i $ H_n({\bf m}).$

\begin{lm}
Елементи $\Delta_i = d^{(\delta)}(u^i),$ $i=2,\ldots, p-1,$ ${u{:=}D(\delta)}$ є нетривіальними,   симетричними інваріантами алгебри  $\overline H_n({\bf m}).$
\end{lm} 
\begin{proof}[Доведення]
Оскільки, $ \mathcal{L}_{0}(u)=0$  то і $ \mathcal{L}_{0}(u^i)=0.$   Тому, з теореми 2 випливає, що в  елементи $d^{(\delta)}(u^i)$  є симетричними інваріантами.  Доведемо  їхню  нетривіальність за індукцією.  Елемент другого степеня  $\Delta_2$   є елементом Казиміра алгебри $ H_n({\bf m}),$(див. \cite{Dz}, \cite{Bk-3})
 і не рівний  нулю. 
Елемент $\Delta_i$ запишемо у вигляді
$$
\begin{array}{l}
\Delta_i = d^{(\delta)}(u^i)=u \Bigl(\sum_{k=1}^{i-2} {i \choose k}u^{k-1} d^{(\delta)}(u^{i-k}) \Bigr)+\Delta_i'=u\Bigl(\sum_{k=1}^{i-2} {i \choose k} u^{k-1}  \Delta_{i-k} \Bigr)+\Delta_i'.
\end{array}
$$
Тут $\Delta_i'= \Delta_i {\bigl | _{u=0}} \in S( H_n({\bf m})).$  Оскільки, за припущенням індукції, всі $\Delta_k,$ $ k<i$ відмінні від нуля і  є алгебраїчно  незалежними  над $\mathbb{K}[u]$, то і $\Delta_i$ відмінний від нуля симетричний інваріант.
\end{proof} \noindent Оскільки, має місце  розклад  $\overline H_n({\bf m})=H_n({\bf m}) \oplus \langle D(\delta) \rangle,$ то кожен елемент  $z$ з $S( \overline H_n({\bf m}))$  однозначно записується у вигляді суми $z=u z_1+z_2,$ де $z_1 \in S( \overline H_n({\bf m})),$ $z_2 \in S(  H_n({\bf m})).$ Якщо ж $z$ -- симетричний інваріант  алгебри $\overline H_n({\bf m}),$ то,  очевидно, $z_2=z{\bigl | _{u=0}}$ буде симетричним інваріантом максимальної підалгебри алгебри $ H_n({\bf m}).$  Таким  чином,  елемент $\Delta_i {\bigl | _{u=0}}  \in  S(  H_n({\bf m}))$  задовольняє умовам теореми 2  і тому він є  генератором симетричного інваріанта  алгебри $ H_n({\bf m}).$ Якщо  ж елемент $\Delta_i {\bigl | _{u=0}},$  рівний  нулю над полем  $ \mathbb{K},$ то розглянемо  відображення  $\varphi :S(L) \to S(L)$ визначене таким  чином $\varphi(a):=\frac{\tau(a)}{p^m}  \mod p,$  де $\tau$  є вкладення  кільця   $S(L)=  \mathbb{K}[L]$  в  кільце $ \mathbb{Z}[L], $  а $m$ --   максимальний   степінь з яким  $p$  входить до розкладу найбільшого  спільного дільника коєфіцієнтів  многочлена $\tau(a) \in  \mathbb{Z}[L].$  Тоді  генератори симетричних  інваріантів будемо  шукати  у вигляді $\varphi (\bigl( \Delta_i {\bigl | _{u=0}})\bigr).$

Нехай $\mathbb{K}^p_n$ -- $p$-центр алгебри $H_n:=H_n(1,1,\ldots,1)$, тобто підалгебра  в $S(H_n)^{H_n}$ породжена  $p$-ми степенями елементів з $S(H_n).$  Покладемо $u_{i,j}=x_1^i x_2^j,$  $(i,j) <(p-1,p-1).$ Неважко переконатися, що для поля характеристики $p=3$ алгебра симетричних інваріантів алгебри $H_2$ рівна $\mathbb{K}^3_2[\Delta_2], $  де \\
$
\Delta_2=d^{(\delta)}(u^2)=(ad(\partial_1))^2 (ad(\partial_2))^2 (u) =2\,{u_{0, \,1}}\,{u_{2, \,1}} + 2\,{u_{1, \,0}}\,{u_{1, \,2}} + {
u_{1, \,1}}^{2} + 2\,{u_{2, \,0}}\,{u_{0, \,2}}.
$

 Наступна теорема описує алгебру симетричних інваріантів алгебри $H_2,$ $ p=5.$
\begin{te} 
 $S(H_2)^{H_2}=\mathbb{K}^5_2[\Delta_2,\Delta_4^* , \Delta_{6}^*],$ де 
$\Delta_i^*= d^{{\phantom{}}(\delta)}\bigl(\varphi (\bigl( \Delta_i {\bigl | _{u=0}}\bigr) \bigr).$
\end{te} 
\begin{proof}[Доведення]
Прямими  обчисленнями  в Maple переконуємося, що   вказані елементи є симетричними  інваріантами, які складаються  відповідно із 12, 78 і 708  доданків.
В  силу громіздскості  випишемо лише  їхні   генератори. 
Генератор симетричного інваріанта $\Delta_2$  рівний $u_{4,4}^2.$  Генератори симетричних інваріантів $\Delta_4^*$ , $\Delta_{6}^*$  мають   відповідно  вигляд 
$$
\begin{array}{l}
 {u_{2, \,3}}^{2}\,{u_{4, \,3}}^{2} + 2\,{u_{4, 
\,2}}\,{u_{4, \,3}}^{2}\,{u_{0, \,4}} + 4\,{u_{0, \,3}}\,{u_{4, 
\,3}}^{3} + 2\,{u_{1, \,2}}\,{u_{4, \,3}}^{2}\,{u_{3, \,4}} + 2\,
{u_{1, \,3}}\,{u_{3, \,3}}\,{u_{4, \,3}}^{2} \\
\mbox{} + 2\,{u_{2, \,1}}\,{u_{4, \,3}}\,{u_{3, \,4}}^{2} + 2\,{u
_{2, \,2}}\,{u_{4, \,3}}^{2}\,{u_{2, \,4}} + 2\,{u_{2, \,2}}\,{u
_{4, \,2}}\,{u_{3, \,4}}^{2} + 2\,{u_{2, \,3}}\,{u_{3, \,3}}^{2}
\,{u_{4, \,3}} \\
\mbox{} + 2\,{u_{2, \,3}}\,{u_{4, \,1}}\,{u_{3, \,4}}^{2} + 2\,{u
_{3, \,1}}\,{u_{3, \,3}}\,{u_{3, \,4}}^{2} + 2\,{u_{3, \,2}}\,{u
_{4, \,3}}^{2}\,{u_{1, \,4}} + 2\,{u_{3, \,2}}\,{u_{3, \,3}}^{2}
\,{u_{3, \,4}} \\
\mbox{} + 2\,{u_{3, \,3}}^{2}\,{u_{4, \,2}}\,{u_{2, \,4}} + 2\,{u
_{4, \,0}}\,{u_{2, \,4}}\,{u_{3, \,4}}^{2} + 2\,{u_{4, \,1}}\,{u
_{4, \,3}}\,{u_{2, \,4}}^{2} + 2\,{u_{4, \,2}}^{2}\,{u_{1, \,4}}
\,{u_{3, \,4}} \\
\mbox{} + 4\,{u_{1, \,3}}\,{u_{4, \,2}}\,{u_{4, \,3}}\,{u_{3, \,4
}} + 4\,{u_{2, \,2}}\,{u_{3, \,3}}\,{u_{4, \,3}}\,{u_{3, \,4}} + 
4\,{u_{2, \,3}}\,{u_{4, \,2}}\,{u_{4, \,3}}\,{u_{2, \,4}} + 4\,{u
_{2, \,3}}\,{u_{3, \,2}}\,{u_{4, \,3}}\,{u_{3, \,4}} \\
\mbox{} + 4\,{u_{2, \,3}}\,{u_{3, \,3}}\,{u_{4, \,2}}\,{u_{3, \,4
}} + 4\,{u_{3, \,1}}\,{u_{4, \,3}}\,{u_{2, \,4}}\,{u_{3, \,4}} + 
4\,{u_{3, \,2}}\,{u_{3, \,3}}\,{u_{4, \,3}}\,{u_{2, \,4}} + 4\,{u
_{3, \,2}}\,{u_{4, \,2}}\,{u_{2, \,4}}\,{u_{3, \,4}} \\
\mbox{} + 4\,{u_{3, \,3}}\,{u_{4, \,1}}\,{u_{2, \,4}}\,{u_{3, \,4
}} + 4\,{u_{3, \,3}}\,{u_{4, \,2}}\,{u_{4, \,3}}\,{u_{1, \,4}} + 
4\,{u_{4, \,1}}\,{u_{4, \,3}}\,{u_{1, \,4}}\,{u_{3, \,4}} + {u_{3
, \,3}}^{4} + 4\,{u_{3, \,0}}\,{u_{3, \,4}}^{3} \\
\mbox{} + {u_{3, \,2}}^{2}\,{u_{3, \,4}}^{2} + {u_{4, \,2}}^{2}\,
{u_{2, \,4}}^{2}\\
\mbox{та } \\
4\,{u_{3, \,3}}\,{u_{4, \,2}}\,{u_{4, \,3}}\,{u
_{3, \,4}}^{3} + 4\,{u_{3, \,3}}\,{u_{4, \,3}}^{3}\,{u_{2, \,4}}
\,{u_{3, \,4}} + 2\,{u_{2, \,3}}\,{u_{4, \,3}}^{3}\,{u_{3, \,4}}
^{2} + 2\,{u_{3, \,2}}\,{u_{4, \,3}}^{2}\,{u_{3, \,4}}^{3} \\
\mbox{} + 3\,{u_{3, \,3}}^{2}\,{u_{4, \,3}}^{2}\,{u_{3, \,4}}^{2}
 + {u_{4, \,1}}\,{u_{4, \,3}}\,{u_{3, \,4}}^{4} + {u_{4, \,3}}^{4
}\,{u_{1, \,4}}\,{u_{3, \,4}} + {u_{4, \,2}}\,{u_{4, \,3}}^{2}\,{
u_{2, \,4}}\,{u_{3, \,4}}^{2} \\
\mbox{} + 3\,{u_{4, \,2}}^{2}\,{u_{3, \,4}}^{4} + 3\,{u_{4, \,3}}
^{4}\,{u_{2, \,4}}^{2}.
\end{array}
$$
Покажемо, що симетричні інваріанти  $\Delta_2,$ $\Delta_4^*,$ $  \Delta_{6}^*$ алгебраїчно незалежні.  Алгебра $H_2$ породжується елементами  вигляду $d^{(\alpha)}(u),$ $\alpha < \delta.$ Визначимо числову адитивну функцію $\lambda$ на $S(H_2)$ таким чином  $\lambda(d^{(\alpha)}(u))=|\alpha|,$ $\lambda(h_1) \lambda(h_2)=\lambda(h_1)+\lambda(h_2),$ $h_1,h_2$ -- мономи з $S(H_2).$ Очевидно, що на   мономах елемента $\Delta_2$ функція $\lambda$ приймає  значення $2(p-1)=8$ для всіх  $i,$  а на  мономах елемента $\Delta_4^*$ i $ \Delta_6^*$ вона приймає значення $4(p-1)=16.$  Прямою  перевіркою знаходимо, що  елементи $\Delta_2^2$  i $\Delta_4^*$ непропорційні. Елемент $6$-го степеня $\Delta_6^*$  може  належати алгебрі $k[\Delta_2,\Delta_4^*]$ лише  у випадку, коли  він є лінійною комбінацією мономів $\Delta_2^3$ i $\Delta_2 \Delta_4^*.$ Але $\lambda(\Delta_2^3)=6 (p-1)=24$  i $\lambda(\Delta_2 \Delta_4^*)=32.$  Оскільки,   маємо ${\lambda( \Delta_6^*)=16,}$ то $\Delta_6^*$   не належить до $k[\Delta_2,\Delta_4^*].$  Таким  чином  вказані три  симетричних  інваріанти є алгебраїчно  незалежними.

Я.С. Крилюком в \cite{Kr} пораховані індекси алгебр Лі картанівського типу, тобто мінімальну кількість нетривіальних  породжуючих елементів центра універсальної огортуючої алгебри. Зокрема  для алгебри $H_2$ індекс рівний $p-2=3.$ Тому знайдена система  із трьох   симетричних  інваріантів $\Delta_2,$ $\Delta_4^*,$ $  \Delta_{6}^*$ породжує нaд $p$-центром алгебру симетричних інваріантів алгебри $H_2.$
\end{proof}

Для  випадку  довільної характеристики поля $p,$ висловимо  наступне припущення

\noindent
{\bf  Гіпотеза.}  {\it $S(H_2)^{H_2}=\mathbb{K}^p_2[\Delta_2,\Delta_4^*,\ldots \Delta_{2(p-2)}^*] .$
}

\end{document}